\documentclass[12pt]{amsart}

\usepackage{amsthm, amsmath, amssymb}
\usepackage{eucal}
\usepackage{graphicx}
\usepackage{color}

\title{Multi-phase Quadrature domains and a related minimization problem}

\author{Avetik Arakelyan}
\address{Institute of Mathematics, NAS of Armenia, Yerevan 0019, Armenia}
\email{arakelyanavetik@gmail.com}

\author{Henrik Shahgholian}
\address{Department of Mathematics, Royal Institute of Technology,
  10044 Stockholm, Sweden}
\email{henriksh@math.kth.se}

\thanks{A. Arakelyan was supported by State Committee of Science MES RA, in frame of the research project No.  SCS 13YR-1A0038}
\thanks{H. Shahgholian was partially supported by the Swedish Research
Council.}
\thanks{Authors are grateful to the referee's careful reading and commenting the first (submitted) version, 
 that resulted in more exact and rigorous statement.}

\newtheorem{theorem}{Theorem}

\newtheorem{lemma}{Lemma}
\newtheorem{proposition}{Proposition}
\newtheorem{corollary}{Corollary}
\theoremstyle{definition}
\newtheorem{definition}{Definition}

\newtheorem{remark}{Remark}

\newcommand{\R}{\mathbb{R}}

 \newcommand{\dist}{\operatorname{dist}}

\begin{document}

\begin{abstract}
In this paper we introduce the multi-phase version of the so-called Quadrature Domains (QD), which refers
to a generalized type of mean value property for harmonic functions. The well-established and developed theory
of one-phase QD was recently   generalized to a two-phase version, by one of the current authors (in collaboration). Here we introduce the concept of
the multi-phase version of the problem, and prove existence as well as several properties of such solutions.
In particular, we discuss possibilities of multi-junction points.
\end{abstract}

\maketitle

\section{Introduction }

A Quadrature domain is a domain that admits a quadrature identity with respect to
a given measure, and a class of functions. The most simple example is any ball in $\R^N$ ($N \geq 2$), which admits  the mean-value property for harmonic functions.
Here the measure is  an appropriate  constant multiple of  Dirac mass at the center of the ball
$$
\int_{B_r(x^0)} h(x) \ dx = < h, c_{N,r}\delta_{x^0}> =c_{N,r}h(x^0)
\qquad \forall \ h \in HL^1(B_r(x^0)),
$$
where $c_{N,r} $ is the volume of the ball $B_r$, and  $HL^1(B_r(x^0))$ denotes the class of integrable harmonic functions over $B_r(x^0):=\{|x-x^0|<r\}$.

More generally, suppose that we are given  a finite, non-negative  measure  $\mu$ with compact support in $\mathbb{R}^N$,  and we want to   find a domain $\Omega \supset \hbox{supp} (\mu)$ such that
\begin{equation}\label{def_QD}
\int_\Omega hdx =(\geq)\int h d\mu,
\end{equation}
for all harmonic (subharmonic)
integrable functions $h$ over the set $\Omega.$ There is vast literature on the topic, which has many connection to diverse areas of analysis. We refer to \cite{gustafsson2005quadrature} and references therein for further background, as well as its connection to other problems.

The concept of two-phase quadrature domain  was introduced  by one of the current authors  (with collaborators) in \cite{MR2754967}, and was further developed  in  \cite{shahgholian2013harmonic}, \cite{gardiner2012two}, \cite{MR2793570}. The interested reader is also referred to  the survey \cite{gardiner2014quadrature},  for the recent advances in quadrature domain theory.

\begin{definition}\label{def:two-phase} (Two-phase QD)
Suppose we are given  constants $\lambda_\pm >0 $,  bounded nonnegative  measures $\mu_\pm,$ and disjoint domains $\Omega_\pm$ such that $supp(\mu_\pm)\subset  \Omega_\pm.$ If for every integrable  harmonic  function  $h$  on $\Omega_+\cup \Omega_-$, that  also has   continuous extension  to $\partial \Omega_+ \cap  \partial \Omega_-$,
  the following integral  identity holds:\footnote{For an inequality as in  \eqref{def_QD}, one needs to assume that 
  $\pm h \chi_{\Omega_\pm}$  is subharmonic.}
\begin{equation}\label{qd-integral}
\int_{\Omega_+} \lambda_+ hdx-\int_{\Omega_-} \lambda_- hdx=(\geq)\int h d(\mu_+-\mu_-),
\end{equation}
then we call $\Omega=\Omega_+\cup \Omega_-$ a \emph{Two-phase quadrature domain} with respect to $\mu_\pm$, and $\lambda_\pm$.
\end{definition}

It is more convenient to    reformulate   the problem in terms of partial differential equation, by using the fundamental solution as the integrand  $h$ above  that results in  the following PDE formulation:
The integral identity  \eqref{qd-integral} is equivalent with the existence of solutions to the following problem (see \cite{MR2754967,gardiner2012two}):

\begin{align}\label{PDE-form}
\begin{cases}
\Delta u= \left( \lambda_+  \chi_{\Omega_+}-  \lambda_- \chi_{\Omega_-} \right) - \left(\mu_+-\mu_-\right)\;\;&\mbox{in}\;\; \mathbb{R}^N,\\
u=0\;\;& \mbox{in}\;\;\mathbb{R}^N\setminus \Omega,
\end{cases}
\end{align}
 where $\chi_A$ stands for the characteristic function of the set $A.$

Examples of non-trivial  two-phase QD are not easily found. One example can be given by solving the one-phase version of the PDE in \eqref{PDE-form}, for a measure $\mu_+$, inside a half-space $\Pi_{\nu, a}=\{ x\cdot \nu > a\}$, with zero boundary values on the plane $\partial  \Pi_{\nu, a}$, and then consider an odd reflection of this solution, which will then produce a two-phase QD for $\mu_+$, and $\mu_-$, where the latter is an (even) reflection of $\mu_+ $ in the same plane.

Although there are yet many unanswered questions concerning two-phase QD, we shall mostly focus on multi-phase version of this problem, and only in Section 3, discuss some geometric results for the two-phase case.
Also,  from time to time we shall need a few known results from both one- and two-phase problems that we shall invoke when needed.

\section{Multi-phase version} \label{sec:multi-phase}

\subsection{The model equation} \label{sec:model}

We shall now  consider a  generalization of the two-phase version, by allowing several phases. To do so, we have to more carefully look into the structure
of a two-phase QD, across the common free boundary, and also the one-phase part of the free boundary. It seems easier to do this from the PDE formulation \eqref{PDE-form}.

  The setting of the problem in terms of partial differential equation is as follows: Given are $m$ positive measures $\mu_i$ and constants $\lambda_i, (i=1,\dots,m).$ We want to find functions $u_i\geq 0, (i=1,\dots,m),$ with
  $\Omega_i \cap \Omega_j = \emptyset$ ($i\neq j$ and $\Omega_i = \{u_i >0\}$) and such that
  \begin{align}
 \Delta(u_i-u_j)=(\lambda_i\chi_{\Omega_i}-\lambda_j\chi_{\Omega_j})-(\mu_i-\mu_j)\;\;&\mbox{in}\;\;\mathbb{R}^N\setminus \cup_{k\neq i,j}\overline{\Omega}_k.
  \end{align}
This, in other words, means that for each pair $(i,j)$ with $i\neq j$ the function
 $u_i - u_j$ solves a two-phase problem outside the union of the supports of the other functions.

Although there are  numerical  constructions of   the multi-phase quadrature domains  for some particular cases (see \cite{MrezaFarid}), the general existence, uniqueness, and stability of such domains so far has been  untouched.

\subsection{The minimizing functional} \label{sec:functional}

In this section we consider certain minimization problem motivated by the recent work \cite{MR2151234}. In \cite{MR2151234} the authors consider a class of stationary states for reaction-diffusion systems of $k\geq 3$ densities having disjoint supports. This problem is arising in population dynamics and turns out to be very useful in studying the problems in potential theory.

We start with the definition of  the minimization sets $K$ and $S.$ Define
\[
K=\{(u_1,u_2,\dots\,u_m)\in (W_0^{1,2}(\mathbb{R}^N))^m,\;\mbox{s.t.}\;\; u_i\geq 0,\; \mbox{for all}\;\;
i=1,\cdots ,m \},
\]
and
\[
S=\{(u_1,u_2,\dots\,u_m)\in (W_0^{1,2}(\mathbb{R}^N))^m,\;\mbox{s.t.}\;\; u_i\geq 0,\;\mbox{and}\;\; u_i\cdot u_j=0,\; \mbox{for all}\;\; i\neq j  \}.
\]
Obviously  we have $S\subset K$. Next we define
\begin{equation}\label{mainfunctinal}
J(u_1,\dots,u_m)=\sum_{i=1}^m\int_{\mathbb{R}^N}\left(\frac{1}{2}|\nabla u_i|^2 - f_i\cdot u_i\right)dx,
\end{equation}
where each function $f_i$  satisfies the following condition.

\vspace{3mm}

\noindent
{\bf Condition A:} (see \cite{MR1390686}):
\begin{equation}
\begin{cases}
\mbox{For all}\;\; i=1,\cdots , m \;\; \mbox{we have}\;\; f_i\in L^\infty(\mathbb{R}^N),\;\emph{supp}(f_i^{+})\;\; \mbox{is compact};\\
f_i\leq\mbox{const.}<0\;\; \mbox{holds outside a compact set.}
\end{cases}
\end{equation}

\begin{remark}
Note that for studying the existence of minimizers to \eqref{mainfunctinal} one could assume more general conditions on $f_i$ rather than only those  in Condition A. E.g., one may let $f_i^-$ to approach zero with reasonable speed as $x$ tends to infinity, and still have minimizers with compact support.
\end{remark}

We need the following lemma.

\begin{lemma}\label{minoverK}
Let all $f_i(x)$ satisfy $Condition$ $A.$ Then $J(u_1,u_2,\dots,u_m)$ is bounded from below and its unique infimum $(v_1,v_2,\dots,v_m)$ is attained  over the set $K.$ Moreover, each component $v_i$ is a unique minimizer to the functional
\[
j_i(U)=\int_{\mathbb{R}^N}\left(\frac{1}{2}|\nabla U|^2 - f_i\cdot U\right)dx,
\]
over the convex set
\[
K_0=\{U\in W_0^{1,2}(\mathbb{R}^N),\;\;\mbox{s.t.}\;\; U\geq 0 \}.
\]
\end{lemma}
\begin{proof}
Observe that $K$ is a closed set, and $J(u_1,u_2,\dots,u_m)$ is lower semi-continuous.
Assume $||f_i||_{\infty}\leq c_i,$ for all $i=1, \cdots , m.$ According to Condition A  there exist  balls $B_{R_i}(0)$ such that
\begin{align*}
J(u_1,u_2,\dots,u_m)&\geq\sum_{i=1}^m\int_{\mathbb{R}^N}\frac{1}{2}|\nabla u_i|^2 dx - \sum_{i=1}^m\int_{B_{R_i}(0)}f_i\cdot u_i dx\\
&\geq\sum_{i=1}^m\int_{B_{R_i}(0)}\frac{1}{2}|\nabla u_i|^2 dx - \sum_{i=1}^m\int_{B_{R_i}(0)}|f_i|\cdot u_i dx\\
&\geq \frac{1}{2}\sum_{i=1}^m\left(||\nabla u_i||_{L^2(B_{R_i}(0))}^2-2c_ip_i||\nabla u_i||_{L^2(B_{R_i}(0))}\right)\geq \sum_{i=1}^m\frac{-(c_ip_i)^2}{2}.
\end{align*}
Here for each  component $u_i$ we have used  Poincare's  inequality $$||u_i||_{L^2(B_{R_i}(0))}\leq p_i||\nabla u_i||_{L^2(B_{R_i}(0))},$$
where the positive constants $p_i$ depend only on the given balls $B_{R_i}(0).$
Thus $J$ is bounded from below, and by lower semi-continuity  $J$ attains it's infimum in $K.$ It is easy to see that $J(u_1,u_2,\dots,u_m)$ is a strictly convex functional. It is well-known that there exists a unique minimizer  $(v_1,v_2,\dots,v_m)$ in $K,$ when $K$ is\emph{ closed and convex} set.

Now, we turn to the second part of the lemma. Indeed, using the previous argument  we can deduce that
the functional
\begin{equation}\label{func_1prmt}
j_i(U)=\int_{\mathbb{R}^N}\left(\frac{1}{2}|\nabla U|^2 - f_i\cdot U\right)dx,
\end{equation}
is coercive and strictly convex, therefore  has a unique minimizer in the\emph{ closed and convex} set $K_0.$ Thus the only thing we have to prove is that for every $i,$ the corresponding component $v_i$ is a  minimizer to \eqref{func_1prmt}. To this end we fix an index  $i_0,$ and assume  $U_{i_0}\in K_0$ is a unique minimizer to \eqref{func_1prmt}.
Apparently, we have $$(v_1,v_2,\dots,v_{i_0-1},U_{i_0},v_{i_0+1},\dots,v_m)\in K.$$ Thus,
\[
J(v_1,v_2,\dots,v_{i_0-1},U_{i_0},v_{i_0+1},\dots,v_m)\geq J(v_1,v_2,\dots,v_m),
\]
which implies
\[
\int_{\mathbb{R}^N}\left(\frac{1}{2}|\nabla U_{i_0}|^2 - f_{i_0}\cdot U_{i_0}\right)dx\geq
\int_{\mathbb{R}^N}\left(\frac{1}{2}|\nabla v_{i_0}|^2 - f_{i_0}\cdot v_{i_0}\right).
\]
According to the definition of  the functional  \eqref{func_1prmt}, we can write the last inequality in the following form:
\[
j_{i_0}(U_{i_0})\geq j_{i_0}(v_{i_0}).
\]
But we know that $U_{i_0}$ is a unique minimizer to \eqref{func_1prmt}, and $v_{i_0}\in K_0$. Thus, we have 
$v_{i_0}=U_{i_0}$,  and this completes the proof of the lemma.
\end{proof}

Now we are ready to prove the main result of this section.
\begin{theorem}\label{embeddtheorem}
Let $f_i(x)$ satisfy $Condition$ $A.$ Then $J(u_1,u_2,\dots,u_m)$ has at least one minimizer
$(\bar{u}_1,\bar{u}_2,\dots,\bar{u}_m)$ in $S$, and also  all minimizers have compact support. 
Moreover, the following  inclusion of supports  holds: 
For any minimizer $(\bar{u}_1,\bar{u}_2,\dots,\bar{u}_m)$ of $J$  over $S$, and 
(the unique) minimizer $(v_1,v_2,\dots,v_m)$ of $J$ over  $K$, we have 
\begin{equation}\label{support}
supp(\bar{u}_i)\subseteq supp(v_i), \qquad  i=1, \cdots , m.
\end{equation}
\end{theorem}
\begin{proof}
From Lemma \ref{minoverK} we know that the functional $J(u_1,u_2,\dots,u_m)$ is lower semi-continuous, coercive and convex. Since the set $S$ is closed, then the existence of a minimizer follows from  standard arguments of calculus of variations. Note that the minimizer is not necessarily unique. Assume that $(\bar{u}_1,\bar{u}_2,\dots,\bar{u}_m)$ is one of the minimizers.

For simplicity, we make the following notations:
\begin{align*}
\bar{U}&\equiv(\bar{u}_1,\bar{u}_2,\dots,\bar{u}_m),\;\;V\equiv(v_1,v_2,\dots,v_m),\\
\min(\bar{U},V)&\equiv(\min(\bar{u}_1,v_1),\min(\bar{u}_2,v_2),\dots,\min(\bar{u}_m,v_m)),\\
\max(\bar{U},V)&\equiv(\max(\bar{u}_1,v_1),\max(\bar{u}_2,v_2),\dots,\max(\bar{u}_m,v_m)).
\end{align*}
To  see the  ordering  of the supports (equation \eqref{support}) for every $i$ we write the following identities, which are easy to verify:
\[
\min(\bar{u}_i,v_i)+\max(\bar{u}_i,v_i)=\bar{u}_i+v_i,
\]
and
\[
\int_{\mathbb{R}^N}\left(|\nabla\min(\bar{u}_i,v_i)|^2+|\nabla\max(\bar{u}_i,v_i)|^2\right) dx=\int_{\mathbb{R}^N}\left(|\nabla\bar{u}_i|^2+|\nabla v_i|^2\right)dx.
\]
Thus, above identities lead us to the following equation
\[
J(\min(\bar{U},V))+J(\max(\bar{U},V))=J(\bar{U})+J(V).
\]
Since $\bar{U}\in S$ and $V\in K,$ then $\min(\bar{U},V)\in S.$ Therefore $$J(\min(\bar{U},V))\geq J(\bar{U}),$$ which implies
$$
J(\max(\bar{U},V))\leq J(V).
$$
Observe that $\max(\bar{U},V)\in K$ and $V\equiv(v_1,v_2,\dots,v_m)$ is a unique minimizer to $J(u_1,u_2,\dots,u_m)$ in $K.$ Hence,
$$
\max(\bar{U},V)=V,
$$
which is equivalent to
\[
(\max(\bar{u}_1,v_1),\max(\bar{u}_2,v_2),\dots,\max(\bar{u}_m,v_m))=(v_1,v_2,\dots,v_m),
\]
therefore $\max(\bar{u}_i,v_i)=v_i,$ for all $i=1, \cdots , m.$ Thus $\bar{u}_i\leq v_i,$ which leads to
\[
supp(\bar{u}_i)\subseteq supp(v_i)
\]
for all $i=1, \cdots , m.$
Due to Lemma \ref{minoverK}, every component $v_i$ is a minimizer to the functional $j_{i}(U)$ in the set $K_0.$ For this type of functionals and  under more general setting  it has been  proved (see \cite[Theorem 1.4]{MR1390686} ) that all minimizers have support in a fixed compact set. Thus, $supp(v_i)$ is compact, which in turn yields the compactness of $supp(\bar{u}_i),$ for all $i=1, \cdots , m.$  This completes the proof of  the theorem.
\end{proof}

\subsection{Special case $m=2$}
The minimization problem we are considering is regarded as a multi-phase free boundary problem. Thus, for the case $m=2,$ it would be natural to expect that the minimization of functional \eqref{mainfunctinal} over the set $S$  is somehow related to the two-phase version considered in  \cite{MR2754967}. For the readers' convenience  we recall the problem discussed in \cite{MR2754967}. The authors considered the following functional
\[
\tilde{J}_{\mathbb{R}^N}(u)=\int_{\mathbb{R}^N}\left\{\frac{1}{2}|\nabla u|^2-f(x)\max(u,0)-h(x)\min(u,0)\right\}dx,
\]
where $f$ and $-h$ satisfy \emph{Condition A}, and
the minimization is taken over the set $ W_0^{1,2}(\mathbb{R}^N) $.
 Using  convexity and coercivity they   obtain the existence of the minimizer $u$.
 In this regard, the minimization functional for our case reads:
\begin{equation}\label{casem=2}
J(u_1,u_2)=\sum_{i=1}^2\int_{\mathbb{R}^N}\left(\frac{1}{2}|\nabla u_i|^2 - f_i\cdot u_i\right)dx,
\end{equation}
where $f_1=f,$ and $f_2=-h.$ The minimization set for the case $m=2$ will be:
\[
S=\{(u_1,u_2)\in (W_0^{1,2}(\mathbb{R}^N))^2,\;\mbox{s.t.}\;\; u_i\geq 0,\;\mbox{and}\;\; u_1\cdot u_2=0\}.
\]

The connection between these two minimization problems is presented in the next theorem.

\begin{theorem}\label{thrmm=2}
Let   $\tilde{u}$ be a minimizer to the functional $\tilde{J}_{\mathbb{R}^N}(u)$ over the set $W_0^{1,2}(\mathbb{R}^N).$ Then
$(\max(\tilde{u},0),-\min(\tilde{u},0))$ is a minimizer to the functional $J(u_1,u_2)$ over the set $S$. Conversely,   if $(\bar{u}_1,\bar{u}_2)$ is a minimizer to the functional $J(u_1,u_2)$  over the set $S,$ then $\bar{u}_1-\bar{u}_2$ is a minimizer to the functional $\tilde{J}_{\mathbb{R}^N}(u)$ over the set $W_0^{1,2}(\mathbb{R}^N)$.
\end{theorem}
\begin{proof}
It is apparent that $(\max(\tilde{u},0),-\min(\tilde{u},0))\in S.$ We have
\begin{align*}
J(\max(\tilde{u},0),-\min(\tilde{u},0))&=\int_{\mathbb{R}^N}\left\{\frac{1}{2}|\nabla \tilde{u}|^2-f(x)\max(\tilde{u},0)-h(x)\min(\tilde{u},0)\right\}dx\\&=\tilde{J}_{\mathbb{R}^N}(\tilde{u}).
\end{align*}
Now, for every $(w_1,w_2)\in S$ we consider  $w=w_1-w_2\in W_0^{1,2}(\mathbb{R}^N).$ Due to $w_1\cdot w_2=0$ we have $w_1=\max(w_1-w_2,0)$ and
$w_2=-\min(w_1-w_2,0).$ Then,
\begin{align*}
\tilde{J}_{\mathbb{R}^N}(\tilde{u})&\leq \tilde{J}_{\mathbb{R}^N}(w_1-w_2)\\&=\int_{\mathbb{R}^N}\left\{\frac{1}{2}|\nabla w|^2-f(x)\max(w_1-w_2,0)-h(x)\min(w_1-w_2,0)\right\}dx\\&= \int_{\mathbb{R}^N}\left(\frac{|\nabla w_{1}|^2}{2}+\frac{ |\nabla w_{2}|^{2}}{2}- f
(x)w_1+ h(x)w_{2} \right)dx,\\&=J(w_1,w_2).
\end{align*}

Thus,
\[
J(\max(\tilde{u},0),-\min(\tilde{u},0))=\tilde{J}_{\mathbb{R}^N}(\tilde{u})\leq J(w_1,w_2),
\]
for every $(w_1,w_2)\in S.$ Hence, $(\max(\tilde{u},0),-\min(\tilde{u},0))$ is a minimizer to the functional $J(u_1,u_2).$

Suppose $(\bar{u}_1,\bar{u}_2)$ is a minimizer to the functional \eqref{casem=2} subject to $S.$ We take arbitrary $w\in W_0^{1,2}(\mathbb{R}^N),$ it is clear that $(\max(w,0),-\min(w,0))\in S.$
Thus, following the same steps of the above computation, yields
\[
\tilde{J}_{\mathbb{R}^N}(\bar{u}_1-\bar{u}_2)=J(\bar{u}_1,\bar{u}_2)\leq J(\max(w,0),-\min(w,0))=\tilde{J}_{\mathbb{R}^N}(w).
\]
Therefore $(\bar{u}_1-\bar{u}_2)$ is a minimizer of the functional $\tilde{J}_{\mathbb{R}^N}(u).$ This completes the proof of the Theorem.
\end{proof}

\begin{remark}
In \cite[Section $2$]{MR2754967} the authors define $U_+$ and $U_-$ to be  the minimizers of functionals  $J_+(u)$ and $J_-(u),$ respectively. Here, for the case $\Omega=\mathbb{R}^N,$ the corresponding functionals are:
\[
J_+(u)=\int_{\mathbb{R}^N}\left(\frac{1}{2}|\nabla u|^2 - f(x)\cdot \max(u,0)\right)dx,
\]
and
\[
J_-(u)=\int_{\mathbb{R}^N}\left(\frac{1}{2}|\nabla u|^2 - h(x)\cdot \min(u,0)\right)dx.
\]
According to Theorem \ref{thrmm=2} we have  $v_1\equiv U_+$ and $v_2\equiv -U_-,$ where $v_i$ is a unique minimizer to the functional $j_i(U)$ defined in Lemma \ref{minoverK}. Thus, we obviously see that the embedding obtained in Theorem \ref{embeddtheorem} is the same as in \cite[Theorem $2.1$]{MR2754967}.
\end{remark}

\subsection{Local properties of minimizers}

In this section we obtain the local properties of the minimizers $(\bar{u}_1,\bar{u}_2,\dots,\bar{u}_m)$ in $S,$ considered in the Theorem \ref{embeddtheorem}. The  next result shows that for every $i\neq j$ the difference $\bar{u}_i-\bar{u}_j$ locally satisfies the two-phase obstacle equation.

\begin{proposition}\label{locprop}
If $(\bar{u}_1,\bar{u}_2,\dots,\bar{u}_m)$  is a minimizer to the functional \eqref{mainfunctinal} subject to the set $S,$ then the following holds in the sense of distributions:
\begin{equation}\label{system1}
\Delta(\bar{u}_i-\bar{u}_j)=-f_i\chi_{\{\bar{u}_i>0\}}+ f_j\chi_{\{\bar{u}_j>0\}}\;\;\qquad \mbox{in}\;\;\mathbb{R}^N\setminus  \cup_{k\neq i,j}\overline{\Omega}_k,
\end{equation}
where $\Omega_i=\{\bar{u}_i>0\}.$ 
\end{proposition}

\begin{proof}
Let  $a^+=\max(a,0)$ and $a^-=-\min(a,0)$.
Take a non-negative test function $\psi \in C_0^\infty(\mathbb{R}^N),$ such that $supp(\psi)\subset\mathbb{R}^N\setminus\cup_{k\neq i,j}\overline{\Omega}_k.$
 For arbitrary $\varepsilon>0$ we set
$$
Q_\varepsilon=\{x\in\mathbb{R}^N,\;\mbox{s.t.}\;\bar{u}_i-\bar{u}_j\leq\varepsilon\psi\}.
$$
 Then we define a new vector $(z_1,z_2,\dots,z_m)$ as follows:
\[
z_l=\bar{u}_l,\;\; \mbox{if}\;\; l\neq i,j\;\; \mbox{and}\;\; z_i=(\bar{u}_i-\bar{u}_j-\varepsilon\psi)^+, z_j=(\bar{u}_i-\bar{u}_j-\varepsilon\psi)^-.
\]
Thus, we clearly  have $(z_1,z_2,\dots,z_m)\in S.$ Keeping in mind that due to $\bar{u}_i\cdot\bar{u}_j=0$ the following identities hold $\bar{u}_i=(\bar{u}_i-\bar{u}_j)^+$ and $\bar{u}_j=(\bar{u}_i-\bar{u}_j)^-,$ we obtain
\begin{align*}
0&\leq J(z_1,z_2,\dots,z_m)-J(\bar{u}_1,\bar{u}_2,\dots,\bar{u}_m)\\&= \int\limits_{\mathbb{R}^N}\left(\frac{1}{2}(|\nabla(\bar{u}_i-\bar{u}_j-\varepsilon\psi)^+|^2-|\nabla\bar{u}_i|^2)\right)+\int\limits_{\mathbb{R}^N}\left(\frac{1}{2}(|\nabla(\bar{u}_i-\bar{u}_j-\varepsilon\psi)^-|^2-|\nabla\bar{u}_j|^2)\right)\\&+
\int\limits_{\mathbb{R}^N}f_i(\bar{u}_i-(\bar{u}_i-\bar{u}_j-\varepsilon\psi)^+)+\int\limits_{\mathbb{R}^N}f_j(\bar{u}_j-(\bar{u}_i-\bar{u}_j-\varepsilon\psi)^-)\\&=
\int\limits_{Q_\varepsilon^c}\left(\frac{1}{2}(|\nabla(\bar{u}_i-\varepsilon\psi)|^2-|\nabla\bar{u}_i|^2)\right)+\int\limits_{Q_\varepsilon}\left(\frac{1}{2}(|\nabla(\bar{u}_i-\bar{u}_j-\varepsilon\psi)|^2-|\nabla\bar{u}_j|^2-|\nabla\bar{u}_i|^2)\right)\\&+
\int\limits_{\mathbb{R}^N}f_i((\bar{u}_i-\bar{u}_j)^+-(\bar{u}_i-\bar{u}_j-\varepsilon\psi)^+)+\int\limits_{\mathbb{R}^N}f_j((\bar{u}_i-\bar{u}_j)^--(\bar{u}_i-\bar{u}_j-\varepsilon\psi)^-)\\&=-\varepsilon\int\limits_{\mathbb{R}^N\setminus\cup_{k\neq i,j}\overline{\Omega}_k}\nabla(\bar{u}_i-\bar{u}_j)\nabla\psi+\frac{1}{2}\varepsilon^2\int\limits_{\mathbb{R}^N\setminus\cup_{k\neq i,j}\overline{\Omega}_k}|\nabla\psi|^2\\&+
\int_{Q_\varepsilon}\left(\frac{1}{2}(|\nabla(\bar{u}_i-\bar{u}_j)|^2-|\nabla\bar{u}_i|^2-|\nabla\bar{u}_j|^2)\right)\\&+
\varepsilon\int\limits_{\mathbb{R}^N\setminus\cup_{k\neq i,j}\overline{\Omega}_k}f_i\chi_{\{\bar{u}_i>\bar{u}_j\}}\psi-\varepsilon\int\limits_{\mathbb{R}^N\setminus\cup_{k\neq i,j}\overline{\Omega}_k}f_j\chi_{\{\bar{u}_i<\bar{u}_j\}}\psi+o(\varepsilon)\\&\leq -\varepsilon\int\limits_{\mathbb{R}^N\cup_{k\neq i,j}\overline{\Omega}_k}\nabla(\bar{u}_i-\bar{u}_j)\nabla\psi+\frac{1}{2}\varepsilon^2\int\limits_{\mathbb{R}^N\setminus\cup_{k\neq i,j}\overline{\Omega}_k}|\nabla\psi|^2\\&+\varepsilon\int\limits_{\mathbb{R}^N\setminus\cup_{k\neq i,j}\overline{\Omega}_k}f_i\chi_{\{\bar{u}_i>\bar{u}_j\}}\psi-\varepsilon\int\limits_{\mathbb{R}^N\setminus\cup_{k\neq i,j}\overline{\Omega}_k}f_j\chi_{\{\bar{u}_i<\bar{u}_j\}}\psi+o(\varepsilon).
\end{align*}
If we divide both sides by $\varepsilon$ and letting $\varepsilon\to 0$ we obtain
\[0\leq -\int\limits_{\mathbb{R}^N\setminus\cup_{k\neq i,j}\overline{\Omega}_k}\nabla(\bar{u}_i-\bar{u}_j)\nabla\psi+\int\limits_{\mathbb{R}^N\setminus\cup_{k\neq i,j}\overline{\Omega}_k}f_i\chi_{\{\bar{u}_i>\bar{u}_j\}}\psi-\int\limits_{\mathbb{R}^N\setminus\cup_{k\neq i,j}\overline{\Omega}_k}f_j\chi_{\{\bar{u}_i<\bar{u}_j\}}\psi.\]
Thus, we have the following inequality in the sense of distributions:
\begin{equation}\label{twophase1}
-\Delta(\bar{u}_i-\bar{u}_j)\leq f_i\chi_{\{\bar{u}_i>\bar{u}_j\}}-f_j\chi_{\{\bar{u}_i<\bar{u}_j\}}\;\;\mbox{in}\;\; \mathbb{R}^N\setminus\cup_{k\neq i,j}\overline{\Omega}_k.
\end{equation}
Interchanging the roles of $i,j$ gives in the
same way

\begin{equation}\label{twophase2}
-\Delta(\bar{u}_j-\bar{u}_i)\leq f_j\chi_{\{\bar{u}_j>\bar{u}_i\}}-f_i\chi_{\{\bar{u}_j<\bar{u}_i\}}\;\;\mbox{in}\;\; \mathbb{R}^N\setminus\cup_{k\neq i,j}\overline{\Omega}_k.
\end{equation}
In view of \eqref{twophase1} and \eqref{twophase2} we will get the two phase equation in the system \eqref{system1}.
\end{proof}

\section{Qualitative and Geometric properties (two-phase case)}
Minimizers to our functional, in the case of one phase problem,
seem to inherit  geometric features that the data enjoys, see   \cite{MR1390686}. In this section we shall  
apply the so-called moving plane method to the two-phase problem, to 
obtain geometric properties of minimizers, that are inherited from data.
The technique seems to fail to be applied to the multi-phase case, due to the simple fact that the class $S$ is not closed under the operation $\max(\bar u_1,\bar u_2)$. Therefore we only consider the two-phase problem in this section. Consider a  fixed unit vector  $n\in \mathbb{R}^N,$ and for $t\in\mathbb{R} $  set
\[
T_t =T_{t,n}=\{x\cdot n=t\}, \;\; T^-_t=T^-_{t,n}=\{x\cdot n<t\}, \;\;\mbox{and}\;\;T^+_t=T^+_{t,n}=\{x\cdot n>t\}.
\]
For $x\in \mathbb{R}^N$ let $x^t$ be a reflected point with respect to $T_t.$ We also set $\varphi^t(x)\equiv\varphi(x^t),$ for a function $\varphi.$ If $\Omega\subset\mathbb{R}^N $ we define
\[
\Omega_t=\Omega\cap T^+_t \;\; \mbox{and} \;\; \tilde{\Omega}_t=\{ x^t\;\mbox{s.t.}\;\; x\in\Omega_t \}.
\]
The following simple lemma will be used in the next theorem

\begin{lemma}(\cite{fp})\label{simple_ineq}
If $\Phi(t)$ is a nondecreasing function of $t\in\mathbb{R},$ and $h_1,h_2$ are $L^\infty$ functions such that $h_1(x)\leq h_2(x).$ Then the following inequality holds:
$$
\int\left(h_1\Phi(z_1)+h_2\Phi(z_2)\right) dx   \leq  \int\left(h_1\Phi(\min(z_1,z_2))+h_2\Phi(\max(z_1,z_2))\right)dx .
$$
\end{lemma}

\begin{theorem}\label{geometr_prop_m=2}
	Let $f_1$, and $f_2$ satisfy Condition A and moreover that for some unit vector $n\in \mathbb{R}^N,$ and some $t_0\in\mathbb{R} $ we have
	\[
	f_1(x)\leq f^t_1(x),\; f_2(x)\geq f^t_2(x)\;\;\mbox{in}\;\; T^+_t, \]
	for all $t\geq t_0.$ Then for every minimizer $(\bar{u}_1,\bar{u}_2)\in S,$ we have
	\begin{align*}
		\bar{u}_1-\bar{u}_2\leq\bar{u}^t_1-\bar{u}^t_2\;\;\mbox{in}\;\; \Omega_t\;\;&\mbox{for all}\;\; t\geq t_0, \\
		\tilde{\Omega}_t\subset\Omega\;\;&\mbox{for all}\;\; t\geq t_0.\\
	\end{align*}
\end{theorem}

\begin{proof}
	We set
	\[
	v^t=
	\begin{cases}
	\min(\bar{u}_1-\bar{u}_2,\bar{u}^t_1-\bar{u}^t_2),\;\;\mbox{in}\;\; T^+_t,\\
	\max(\bar{u}_1-\bar{u}_2,\bar{u}^t_1-\bar{u}^t_2)\;\;\mbox{in}\;\; T^-_t.
	\end{cases}
	\]
	Let
	\[
	L(\varphi)=\int_{T^+_t} \left(\frac{1}{2}|\nabla \varphi|^2-f_1\varphi^+-f_2\varphi^-\right)dx,
	\]
	and
	\[
	L_t(\varphi)=\int_{T^+_t} \left(\frac{1}{2}|\nabla \varphi|^2-f^t_1\varphi^+-f^t_2\varphi^-\right)dx.
	\]

	Now if we apply Lemma \ref{simple_ineq} one time for $z_1=\bar{u}_1-\bar{u}_2,$ $z_2=\bar{u}_1^t-\bar{u}_2^t,\;$ $h_1=f_1,\;$$h_2=f_1^t$ and $\Phi(t)=\max(t,0),$ and the second time for  $z_1=\bar{u}_1-\bar{u}_2,$ $z_2=\bar{u}_1^t-\bar{u}_2^t,\;$ $h_1=-f_2,\;$ $h_2=-f_2^t$ and $\Phi(t)=\min(t,0)$ we will get
	$$
		\int_{T^+_t} \left( f_1 z_1^+  +  f_1^t z_2^+ \right) dx \leq  
		\int_{T^+_t}  \left( f_1 (\min(z_1,z_2))^+ +   f_1^t (\max(z_1,z_2))^+ \right)dx,
	$$
	and
	$$
    \int_{T^+_t} \left( f_2 z_1^-  +  f_2^t z_2^- \right) dx \leq \int_{T^+_t}  \left( f_2(\min(z_1,z_2))^- +   f_2^t (\max(z_1,z_2))^- \right) dx.
	$$
	Hence,
	\begin{align*}
		\tilde{J}_{\mathbb{R}^N}(v^t)&=L(\min(\bar{u}_1-\bar{u}_2,\bar{u}^t_1-\bar{u}^t_2))+L_t(\max(\bar{u}_1-\bar{u}_2,\bar{u}^t_1-\bar{u}^t_2))\\
		&\leq L(\bar{u}_1-\bar{u}_2)+L_t(\bar{u}^t_1-\bar{u}^t_2)=\tilde{J}_{\mathbb{R}^N}(\bar{u}_1-\bar{u}_2),
	\end{align*}
	for all $t\geq t_0.$ Thus, $v^t$ is  also a minimizer to the functional $\tilde{J}_{\mathbb{R}^N}.$ It is easy to see that the minimizer of  $\tilde{J}_{\mathbb{R}^N}$ must be unique due to its  coercivity and convexity. Hence, $v^t\equiv\bar{u}_1-\bar{u}_2,$ which yields  $$\bar{u}_1-\bar{u}_2\leq\bar{u}^t_1-\bar{u}^t_2;\;\mbox{in}\;\;T^+_t\;\;\mbox{for all}\;\; t\geq t_0,$$
	and $\tilde{\Omega}_t\subset\Omega$ for all $ t\geq t_0.$ This also implies that $n\cdot\nabla(\bar{u}_1-\bar{u}_2)\leq 0$ in $\Omega_{t_0}.$
\end{proof}

\begin{corollary}
	Let  $(\bar{u}_1,\bar{u}_2)\in S,$ and $f_i(x), i=1,2$ be as in Theorem \ref{geometr_prop_m=2}. If we assume  $f_i(x), i=1,2$ are symmetric in $T_{t_0},$ then $\bar{u}_1-\bar{u}_2$ is symmetric in $T_{t_0}.$
\end{corollary}

Another observation can be made when the ingredients have scaling properties, that are inherited by  solutions.
This is reflected in our next result.

\begin{theorem}
	Let  $(\bar{u}_1,\bar{u}_2)\in S,$ is a minimizer to \eqref{casem=2} over  the set $S,$ and  $f_1, f_2$  satisfy Condition A. If, moreover we assume that
	\[
	t^\alpha f_1(x/t)\leq f_1(x),\;t^\alpha f_2(x/t)\geq f_2(x)\;\;\mbox{for all}\;\; t\in (0,1) \;\;\mbox{and}\;\; x\in\mathbb{R}^N,
	\]
	then
	\[t^{\alpha+2}(\bar{u}_1(x/t)-\bar{u}_2(x/t))\leq \bar{u}_1(x)-\bar{u}_2(x),\]
	for every fixed real number $\alpha.$ In particular, for the case $\alpha=-1,$ we get that  $t(\bar{u}_1(x/t)-\bar{u}_2(x/t))\leq \bar{u}_1(x)-\bar{u}_2(x),$ which in turn implies that the set $\{\bar{u}_1(x)>\bar{u}_2(x)\}=\{\bar{u}_1(x)>0\}$ is starshaped with respect to the origin.
\end{theorem}

\begin{proof}
First, observe that if $(\bar{u}_1(x),\bar{u}_2(x))$ is a minimizer to \eqref{casem=2}, then
$$(t^{\alpha+2}\bar{u}_1(x/t),t^{\alpha+2}\bar{u}_2(x/t))$$ is going to be a minimizer to

\begin{equation}\label{homogen_m=2}
J_t(v_1,v_2)=\sum_{i=1}^2\int_{\mathbb{R}^N}\left(\frac{1}{2}|\nabla v_i|^2 - t^{\alpha}f_i(x/t)\cdot v_i\right)dx,
\end{equation}
subject to the set $S.$ Thus, applying Theorem \ref{thrmm=2}  we obtain that  $t^{\alpha+2}(\bar{u}_1(x/t)-\bar{u}_2(x/t))$ is a  minimizer to the following  functional
\[
\tilde{J}^t_{\mathbb{R}^N}(u)=\int_{\mathbb{R}^N}\left\{\frac{1}{2}|\nabla u|^2-t^{\alpha}f_1(x/t)\max(u,0)+t^{\alpha}f_2(x/t)\min(u,0)\right\}dx,
\]
subject to the set $ W_0^{1,2}(\mathbb{R}^N).$ Now, repeating the similar arguments as those  in the proof of Theorem \ref{geometr_prop_m=2} we will get that
\[t^{\alpha+2}(\bar{u}_1(x/t)-\bar{u}_2(x/t))\leq \bar{u}_1(x)-\bar{u}_2(x),\]
for every $t\in(0,1).$ And for the case $\alpha=-1$ this readily shows that  $\{\bar{u}_1(x)>\bar{u}_2(x)\}=\{\bar{u}_1(x)>0\}$ is starshaped with respect to the origin. This completes the proof of Theorem.
\end{proof}

\section{Multi-phase Quadrature Identity} \label{sec:mpqi}
\subsection{Minimization with mollified Radon measures}
Throughout this section we will consider the existence of solutions  to  \eqref{mainfunctinal} for the  case when each $f_i$ is replaced by $\mu_i \ast \psi -\lambda_i$, where $\mu_i\ast \psi $ is a mollified version of $\mu_i$ which is a  positive Radon measure with compact support and $\lambda_i>constant  >0$ are positive $L^\infty$ functions.
Our main concern is to consider the  existence of solutions satisfying the system \eqref{system1}, with $f_i=\mu_i-\lambda_i.$

In the sequel, the following approximation theorem,  will play a crucial role.

\begin{theorem}(\cite{MR2754967},Theorem $3.3$)\label{Shahghth1}
Let $f=\mu-\lambda$ where $\mu$ is a Radon measure with compact support, 
 and let $\varphi_n$ be a sequence of smooth functions with compact support, contained in a a fixed compact set, and  such that $\varphi_n\to \mu$ weakly as measures. Let $u_n$ be the minimizer for the functional
\[
E_n(u)=\int_{\mathbb{R}^N}\left(\frac{1}{2}|\nabla u|^2 - (\varphi_n(x)-\lambda(x))u\right)dx,
\]
over the set $K_0.$
Then there exists a compact set $F$ such that $supp(u_n)\subset F$ for all $n.$
\end{theorem}

The main result of this section reads as follows.

\begin{theorem}\label{Radon_sol}
	Given are $m$ positive Radon measures $\mu_i$  with compact support and bounded functions $\lambda_i$ such that $\lambda_i > const.>0.$ Then  there exists at least one $(u^*_1,u^*_2,\dots,u^*_m)\in S$ solution of the system  \eqref{system1}, with $f_i=\mu_i-\lambda_i,\;\forall i=1, \cdots , m.$ 
\end{theorem}

\begin{proof}
For every measure $\mu_i$ we take the sequence $\varphi^n_i$ of smooth functions with compact support such that $\varphi^n_i\to \mu_i$ as defined in Theorem \ref{Shahghth1}. It is clear that for all $n$ and every $i$ the differences $f^n_i=\varphi^n_i-\lambda_i$ will satisfy the Condition A, and $f^n_i\to (\mu_i-\lambda_i)$ in the sense defined in Theorem \ref{Shahghth1}. We consider the following functional:
\begin{equation}\label{mainfunctinal2}
J_n(u_1,\dots,u_m)=\sum_{i=1}^m\int_{\mathbb{R}^N}\left(\frac{1}{2}|\nabla u_i|^2 - (\varphi^n_i-\lambda_i)\cdot u_i\right)dx.
\end{equation}

Due to Theorem \ref{embeddtheorem} the functional \eqref{mainfunctinal2} has at least one minimizer subject to $S,$ which we denote by $(\bar{u}^n_1,\bar{u}^n_2,\dots,\bar{u}^n_m).$ According to Proposition \ref{locprop} the minimizer satisfies the system \eqref{system1}, with right-hand side  $f^n_i,$ in the distributional sense.
 Theorem \ref{embeddtheorem} also implies that for all $n$ we have the following embedding
\[
supp(\bar{u}^n_i)\subseteq supp(v^n_i),
\]
for all $i=1, \cdots , m$, where $(v^n_1,v^n_2,\dots,v^n_m)$ is the unique minimizer to $J_n(u_1,u_2,\dots,u_m)$ over the set $K.$
Now, in view of Lemma \ref{minoverK} we have that every component $v^n_i$ is a unique minimizer to the following functional:
\begin{equation}\label{minoverK2}
j^n_i(U)=\int_{\mathbb{R}^N}\left(\frac{1}{2}|\nabla U|^2 - (\varphi^n_i-\lambda_i)U\right)dx,
\end{equation}
over the set $K_0.$ Thus,  for every $i$ due to Theorem \ref{Shahghth1} there exists a compact set $Q_i$ such that for all $n$ we have
\[
supp(\bar{u}^n_i)\subseteq supp(v^n_i)\subset Q_i.
\]
Hence, the support of the minimizers $(\bar{u}^n_1,\bar{u}^n_2,\dots,\bar{u}^n_m)$  remain in a compact set $\cup_{i=1}^m Q_i$ in the limit. This implies that there exists a subsequence which is $weak^*$-convergent as distributions to the limit  $(u^*_1,u^*_2,\dots,u^*_m),$ which will clearly satisfy the system \eqref{system1}, in the distributional sense, with right-hand side $\mu_i-\lambda_i,$ for all $i=1, \cdots , m.$ This completes the proof of Theorem.
\end{proof}

Another interesting approach of solution existence for  the system \eqref{system1}, when $f_i=\mu_i-\lambda_i,$ can be performed for certain sufficiently concentrated (in the sense  defined in \cite{MR1390686}) measures $\mu_i.$ To this aim, we are going to use   Lemma $4.6$ in \cite{MR1390686}.
  Thus, as in the previous Theorem \ref{Shahghth1}, we consider the minimization problem for the case  $f_i=\mu_i\ast\psi-\lambda_i$ for every fixed $i.$ Here, the mollifier $\psi$ is chosen such that $0\leq\psi\in L^\infty(\mathbb{R}^N)$ be a non-increasing, radially symmetric function satisfying $\int\psi dx=1.$

 The result for the concentrated measures considered in \cite[Theorem $4.7$]{MR1390686} will be the following theorem below:

\begin{theorem}
Let the positive measures $\mu_i$ be  sufficiently concentrated in some balls $B(x_i,R_i)$  such that
\begin{equation}
	\begin{cases}
		\mu_i(B^c_{R_i})=0,\\
		\mu_i(B_{R_i})>(b_i+\frac{Nc_i}{3{R_i}})6^N|B_{R_i}|,
	\end{cases}
\end{equation}
where $b_i,c_i\geq 0$ are constants with $b_i+c_i>0,$ for all $i=1,2,\dots,m.$ Then there exists $(\tilde{u}_1,\tilde{u}_2,\dots,\tilde{u}_m)\in S$ solution of the system \eqref{system1}, with $f_i=\mu_i-\lambda_i,\;\forall i=1, \cdots , m,$ where $\lambda_i>const.>0$ are $L^\infty$ functions.
\end{theorem}

\begin{proof}
Assume that there exists constants $l_i>0$ such that $\lambda_i(x)>l_i>0,$ and the measures $\mu_i$ are sufficiently concentrated in balls $B(x_i,R_i)$.  Theorem $4.7$ in \cite{MR1390686} provides the existence of 
$(w_1,w_2,\dots,w_m)\in K$ which satisfy (in the sense of distributions)
\begin{equation}
\begin{cases}
\Delta w_i=\lambda_i-\mu_i\;\;\mbox{in}\;\;\{w_i>0\}\\
w_i=|\nabla w_i|=0\;\;\mbox{on}\;\;\partial\{w_i>0\}\\
supp(\mu_i)\subset\{w_i>0\},
\end{cases}
\end{equation}
for every $i=1, \cdots , m$. According to Lemma $4.6$ in \cite{MR1390686} every component $w_i$ can be viewed as a minimizer of the following functional
\[
E^\psi_i(u)=\int_{\mathbb{R}^N}\left(\frac{1}{2}|\nabla u|^2 - ((\mu_i\ast\psi)(x)-\lambda_i(x))u\right)dx,
\]
subject to the set $K_0.$ Here $\psi$ is defined as in \cite[Lemma $4.6$]{MR1390686}. We can always choose the radially symmetric function $\psi=\psi_r$ such that $(\mu_i\ast\psi_r)(x)-\lambda_i(x)$ satisfies the Condition A. Thus, the minimizer to the functional $E^\psi_i(u)$ will have a compact support. After applying Theorem \ref{embeddtheorem} we obtain $supp(u_i^{r})\subseteq supp(w_i),$ where $(u_1^{r},u_2^{r},\dots,u_m^{r})$ is a minimizer of
\[
E^\psi(u_1,u_2,\dots,u_m)=\sum_{i=1}^m\int_{\mathbb{R}^N}\left(\frac{1}{2}|\nabla u_i|^2 - ((\mu_i\ast\psi)(x)-\lambda_i(x))u_i\right)dx,
\]
subject to the set $S,$ for $\psi=\psi_r.$ Thus, the vector $(u_1^{r},u_2^{r},\dots,u_m^{r})$ remain in a compact set. Again as in Theorem \ref{Radon_sol} one can pass to the limit $r\to 0,$ and since the function $\psi_r$ has a compact support in $B_r(0),$ then apparently we will obtain $(\mu_i\ast\psi_r)(x)\to \mu_i,$ and $u_i^{r}\to\tilde{u}_i.$ This completes the proof of Theorem.
\end{proof}

\begin{remark}
Observe that although for every $i$ we have $supp(\mu_i)\subset\{w_i>0\}$. This may easily fail for the multi-phase case, i.e. we in general can not claim  $supp(\mu_i)\subset supp(\tilde{u}_i),$ where $(\tilde{u}_1,\tilde{u}_2,\dots,\tilde{u}_m)\in S$ is the solution of the system \eqref{system1}. 
\end{remark}

\subsection{Multi-phase QD } \label{sec:mpqd}

In this section we  discuss the concept of Balayage and Quadrature Domains for the general (multi-phase) case.
The one phase problem has been well studied in the literature and we refer the reader to the following works \cite{gustafsson1990quadrature, gustafsson2005quadrature, sakai1982quadrature, sakai1983applications}. The concept of two phase version of the so-called Quadrature domain has been defined and studied recently, in the works \cite{MR2754967,gardiner2012two}.

The key point is that the measures have to be concentrated enough and also  in balance. Indeed, if one of the measures $\mu_i$, say, has a very high density on its support, but not the others, then the support of the corresponding $u_i$ will have the possibility of  covering  the support of $\mu_j$, for $j \neq i$.  This naturally makes 
it impossible to find an $m$-QD for our measures. 
 Finding right conditions for this balance is a question to be answered in the future. Here we will illustrate this for measures that satisfy Sakai's condition. First of all, we take $\lambda_i>0$ to be constant for all $i=1,\dots, m.$  The main difficulty is to provide  conditions which lead to the existence of solution of the system \eqref{system1} with property $supp(\mu_i)\subset \{\overline{u}_i>0\},$ for all $i=1, \cdots , m.$ The latter property  implies the following condition $\mu_i\equiv\mu_i\chi_{\{\bar{u}_i>0\}},$ for all $i=1, \cdots , m.$ Thus,  system \eqref{system1}, for $f_i=\mu_i-\lambda_i,$  and $\Omega_i=\{u_i >0 \}$  can be rewritten as follows:

\begin{equation}\label{system2}
\Delta(\bar{u}_i-\bar{u}_j)=(\lambda_i\chi_{\Omega_i} -\lambda_j\chi_{ \Omega_j }  )-(\mu_i-\mu_j)\;\;\mbox{in}\;\;\mathbb{R}^N\setminus\cup_{k\neq i,j}\overline{\Omega}_k.
\end{equation} 
For every $i\neq j$ let   $h \in HL^1 (\Omega_i \cup \Omega_j) $, which  has also    continuous extension  to 
$\Omega_i \cup \Omega_j \cup (\partial \Omega_i \cap \partial \Omega_j) \cup (\cup_{k \ne i,j} \partial \Omega_k)$, 
and $h=0$ on $ \cup_{k\neq i, j} \partial \Omega_k.$ 
Next we formally write (leaving the verification to the reader\footnote{To verify this one needs to show 
$\int_{\Omega_i \cup \Omega_j }h \Delta(\bar{u}_i -\bar{u}_j)=0$. This can be done using a sequence of smooth $C^\infty$ functions $\omega_k$ with compact support in $ interior(\overline{\Omega_i \cup\Omega_j})$ and further properties as described in \cite{sakai1983applications}, the proof of Proposition 4. It should be remarked that for the boundary part $(\partial \Omega_i \cup \partial \Omega_j) \cap \partial \Omega_k$ where $k\neq i, j$ one has to use the assumption that $h=0$ there.
})

\[
\int_{\Omega_i \cup \Omega_j }h(\lambda_i\chi_{\Omega_i }-\lambda_j\chi_{\Omega_j })
=\int_{\Omega_i \cup \Omega_j }h(\Delta(\bar{u}_i-\bar{u}_j)+(\mu_i-\mu_j)),
\]
which (using the properties of $h$ mentioned above) leads to
\[
\int_{\Omega_i }h\lambda_i-\int_{\Omega_j }h\lambda_j=\int h(\mu_i-\mu_j).
\]
It is easy to see that the standard  mollifier technique (see \cite{MR1390686}) will also work in this case, and we may replace the measures with smooth functions, with support close to the support of measures.

\begin{definition}[Multi-phase Quadrature domain]
Suppose we are given $m$ bounded positive measures $\mu_i,$ and disjoint domains $\Omega_i$ such that $supp(\mu_i)\subset \Omega_i.$ For each $i\neq j$ let   $h \in HL^1 (\Omega_i \cup \Omega_j) $, $h$ is continuous across $\partial \Omega_i \cap \partial\Omega_j$,  
and $h=0$ on $ \cup_{k\neq i, j} \partial \Omega_k.$ If for each $i \neq j$ the above class of  harmonic functions  admit the  following QI
\begin{equation}\label{Q-identity}
\int_{\Omega_i}h\lambda_idx-\int_{\Omega_j}h\lambda_jdx= \int h d(\mu_i-\mu_j),
\end{equation}
then we call $\Omega=\{\Omega_i\}_{i=1}^m$ an  {\it  $m$-phase QD} with respect to the measure $\{\mu_i\}_{i=1}^m$, and the positive constants $\{\lambda_i\}_{i=1}^m$. (In general $\lambda_i$ can be taken to be strictly positive 
functions.)

If we reduce the  test class $h$ to be subharmonic in $\Omega_i$ and super-harmonic in $\Omega_j$ (due to negative sign in front of the integral) then  the equality  in \eqref{Q-identity} is replaced with an inequality $(\geq) $.

Observe also if $m=1$, then we may take $\lambda_2= \dots  =\lambda_m=0$, $\Omega_i=\emptyset$ for $i=2, \dots , m$,  and we have the definition of a one-phase quadrature domain
\begin{equation}\label{1-Q-identity}
\int_{\Omega_1}h\lambda_1 dx  = \int h d\mu_1.
\end{equation}
\end{definition}

\begin{definition}[Sakai's concentration condition]
We say that the measures $\mu_i$ satisfy \emph{Sakai's concentration condition} if for every $i$ and $x\in supp(\mu_i)$
\[
\underset{r\to 0^+}{\;\limsup}{\frac{\mu_i(B_r(x))}{|B_r|}}\geq 2^N\lambda_i,
\]
where $\lambda_i$ are positive constants for $i=1, \cdots , m$.
\end{definition}

\begin{theorem}\label{existenceQD}
Let $\mu_i$ be  given Radon measures with compact supports, and $\lambda_i$ are positive constants that satisfy Sakai's concentration condition. Suppose that for each $\mu_i$ the corresponding  one-phase quadrature domain 
  $Q_i$ (see \eqref{1-Q-identity})  is such that
\begin{equation}\label{disjointQD}
\overline{Q}_i\cap supp(\mu_j)=\emptyset,\;\; \mbox{for every}\;\; i\neq j.
\end{equation}
 Then, we have a solution to our multi-phase free boundary problem \eqref{system2}.
\end{theorem}

\begin{proof}

We  consider  mollifiers  $\mu_i \ast \psi $, of the measures   $\mu_i$  ($i=1, \cdots , m$) and minimize the functional  \eqref{mainfunctinal} for $f_i = \mu_i \ast \psi -\lambda_i $. Since $supp(\mu_i \ast \psi ) $ is a subset of a $\epsilon$-neighborhood of $supp(\mu_i)$, it suffices (by taking $\epsilon $ arbitrary small) to show the theorem for smooth $\mu_i$. We thus from now on assume $\mu_i$ is smooth enough such that a minimizer 
 $(\bar{u}_1,\bar{u}_2,\dots,\bar{u}_m)$ belongs to  $ S$.

Next, if we prove that $supp(\mu_i)\subset supp(\bar{u}_i),$ then obviously $(\bar{u}_1,\bar{u}_2,\dots,\bar{u}_m)$ will solve the multi-phase free boundary problem \eqref{system2}. To this end, first observe that due to Sakai condition we have the following embedding $supp(\mu_i)\subset Q_i,$ for all $i=1, \cdots , m.$ We argue by contradiction.
Assume that there exists a number $i_0$ such that $supp(\mu_{i_0})\setminus supp(\bar{u}_{i_0})\neq\emptyset.$ Then, according to condition \eqref{disjointQD} and Theorem \ref{embeddtheorem}, there exists a point $z_0\in supp(\mu_{i_0})\setminus supp(\bar{u}_{i_0}),$ such that $\dist(z_0,\tilde{\Omega})>0,$ where $\tilde{\Omega}=\cup_{i=1}^m supp(\overline{u}_i).$ Thus, one can easily take a ball $B_R(z_0)$ such that $\overline{B_R(z_0)}\cap\tilde{\Omega}=\emptyset.$

Let $\underset{i}{\min}\;\lambda_i\geq\varepsilon>0,$ and $\underset{i}{\sup}\; |\mu_i|\leq M<\infty$ such that $\varepsilon<M.$ We define $r=\left(\frac{\varepsilon}{M}\right)^{1/N}\cdot R$ and consider the following measure, which  satisfies Sakai's condition with respect to $\lambda_{i_0}$
\[
\nu_{i_0}\equiv\mu_{i_0}\cdot\chi_{B_{r}(z_0)}.
\]
We have
\begin{equation}\label{ineq}
\nu_{i_0}-\lambda_{i_0}\leq M\cdot\chi_{B_r(z_0)}-\lambda_{i_0}\leq M\cdot\chi_{B_r(z_0)}-\varepsilon.
\end{equation}

Define
\[
L_{M,\varepsilon}(U)=\int_{\mathbb{R}^N}\left(\frac{1}{2}|\nabla U|^2 - (M\cdot\chi_{B_r(z_0)}-\varepsilon)U\right)dx.
\]
According to Lemma $1.2$ in \cite{MR1390686} every minimizer of $L_{M,\varepsilon}(u)$ over the set $K_0$ is radially symmetric, radially non-increasing and vanishes outside a compact set. Moreover, the minimizer is unique (in our case $c=0$) and its support is a ball centered at $z_0$ and with radius $\sigma =\left(\frac{M}{\varepsilon}\right)^{1/N}\cdot r=R$ (see Example $1.5$ in \cite{MR1390686}).  The proof of this result relies on the so-called symmetric decreasing rearrangement technique, and we refer for its background to the book \cite{mossino1984inegalites}. Using the same arguments as in the proof of  Theorem $1.4$ in \cite{MR1390686}, and inequality \eqref{ineq}, one can easily conclude that   $supp(v_{i_0})\subset \overline{B_R(z_0)}.$ Where $v_{i_0}$  is a unique minimizer to the functional
\[
J_{\nu_{i_0},\lambda_{i_0}}(U)=\int_{\mathbb{R}^N}\left(\frac{1}{2}|\nabla U|^2 - (\nu_{i_0}-\lambda_{i_0})U\right)dx,
\]
over the set $K_0.$ Observe that the set  $\{v_{i_0}>0\}$ is a one-phase Quadrature domain with respect to the measure $\nu_{i_0}$ and constant $\lambda_{i_0}$ (see Section $4$ in \cite{MR1390686}). Since $\nu_{i_0}$ satisfies Sakai's concentration, then the set  $\{v_{i_0}>0\}$ is not empty, namely $v_{i_0}\neq 0.$ Moreover, due to $J_{\nu_{i_0},\lambda_{i_0}}(0)=0,$ we clearly have $J_{\nu_{i_0},\lambda_{i_0}}(v_{i_0})<0,$ which will be used later on.

Thus, $\overline{B_R(z_0)}\cap\tilde{\Omega}=\emptyset$ implies $supp(v_{i_0})\cap\tilde{\Omega}=\emptyset,$ and therefore
\[
W\equiv\left(\bar{u}_1,\bar{u}_2,\dots,\bar{u}_{i_0-1},\bar{u}_{i_0}+v_{i_0},\bar{u}_{i_0+1},\dots,\bar{u}_m\right)\in S.
\]
Hence,
\begin{align*}
J(W)&=
\sum_{i\neq i_0}\int_{\mathbb{R}^N}\left(\frac{1}{2}|\nabla \bar{u}_i|^2 - (\mu_i-\lambda_i)\bar{u}_i\right)dx\\&+\int_{\mathbb{R}^N}\left(\frac{1}{2}|\nabla (\bar{u}_{i_0}+v_{i_0})|^2 - (\mu_{i_0}-\lambda_{i_0})(\bar{u}_{i_0}+v_{i_0})\right)dx\\&=J\left(\bar{u}_1,\bar{u}_2,\dots,\bar{u}_m\right)+\int_{\mathbb{R}^N}\left(\frac{1}{2}|\nabla v_{i_0}|^2 - (\mu_{i_0}-\lambda_{i_0})v_{i_0}\right)dx\\&\leq J\left(\bar{u}_1,\bar{u}_2,\dots,\bar{u}_m\right)+\int_{\mathbb{R}^N}\left(\frac{1}{2}|\nabla v_{i_0}|^2 - (\nu_{i_0}-\lambda_{i_0})v_{i_0}\right)dx\\&< J\left(\bar{u}_1,\bar{u}_2,\dots,\bar{u}_m\right),
\end{align*}
which contradicts the minimality of $(\bar{u}_1,\bar{u}_2,\dots,\bar{u}_m).$ Thus $ supp(\mu_{i_0})\setminus supp(\bar{u}_{i_0})=\emptyset,$ and this implies  $supp(\mu_{i_0})\subset supp(\bar{u}_{i_0}).$ This completes the proof of Theorem.
\end{proof}

\section{Analysis of junction points}
In this section we shall discuss  junction points, where several phases meet. We shall keep the discussion partially at an informal and heuristic level, as the analysis needed to support our argument could be quite  complicated and outside the scope of this paper. However, our main results Theorems \ref{thm-three-phase} and \ref{four-junction}
are proven with complete mathematical  rigor. 

 It is also obvious that only  junction points away from the supports of the measures 
are interesting, and subject for study. Nevertheless, we shall also give examples when the support of the measures hit a junction point  (see Figure  \ref{QDfour}).

Let us first remark that  if the  measure doesn't hit the junction point,  then there is a quadratic non-degeneracy for solutions to our problem, due to $\Delta u_i = \lambda_i$ in the set $\Omega_i$, i.e.  one has
$$
\sup_{B_r(x^0)}  u_i (x) \geq \frac{\lambda_i}{2n} r^2 +  u_i (x^0), \qquad \forall \ x^0 \in \Omega_i
$$
by a simple application of the maximum principle to the function $ u_i (x)-  u_i (x^0) - \frac{\lambda_i}{2n}|x-x^0|^2$ in the domain $\Omega_i \cap B_r(x^0) $
(see \cite{caffarelli1977regularity}). Letting $x^0$ tend to a junction point we have the quadratic  non-degeneracy
\begin{equation}\label{non-deg}
\sup_{B_r(x^0)}  u_i (x) \geq \frac{\lambda_i}{2n} r^2 , \qquad \forall \ x^0 \in \partial  \Omega_i .
\end{equation}

This non-degeneracy, along with a barrier argument  implies that the set $\Omega_i \cap B_r(x^0), $ where $x^0$ is a junction point, cannot be to narrow in the sense that it cannot be confined within a cone of aperture  less than $\pi/2$ (see footnote below.\footnote{The reader familiar with the topic of minimal partition and segregation problem, will see a qualitative difference of the analysis of junctions for our problem and minimal partition, in the sense of the "degree" of the junction.})
 Indeed, in such cones we have nonnegative (homogeneous)  harmonic functions $h$ that decay faster than quadratic, we can always make $ u_i \leq Ch$ (for some large $C$).  But then the quadratic non-degeneracy is violated, and we have a contradiction.

The above argument seems also to work for domains $\Omega_i$, where $\Omega_i \cap \partial B_r(x^0) $
is small enough such that the corresponding eigenvalue problem on the sphere $\partial B_r(x^0)$ admits the first
eigenvalue strictly larger than 2, for all $r \leq r_0$ (for some $r_0 >0)$. This then allows for a barrier and a contradiction argument through maximum principle.

The, partially heuristic, analysis above show that at a junction point (which is outside the support of the measure) each domain $\Omega_i$ cannot be confined in a cone of small size (in two dimensions the cone cannot have an angle  smaller than $\pi/2$). This suggests that, formally, 
we cannot have junction points (in two-dimensions) where the support of five solutions meet,  provided the junction point doesn't hit the measure. It is not straightforward   whether a junction point of four solution can exists or not, as this seems to be a border line, see Section \ref{sec-4-junction}.
 The  junctions point with three phases are possible, which is  proven below, in Section \ref{three-j}.

\begin{remark}
Note that when the junction point  hits the measure we will have  $\Delta u_i = \lambda_i-\mu_i$ close to the junction point, instead of  $\Delta u_i = \lambda_i.$ Thus, the non-degeneracy argument fails, and it might be possible  that the set  $\Omega_i \cap B_r(x^0), $ where $x^0$ is a junction point,  is narrow enough in the sense that it can be confine within a cone of aperture less than $\pi/2.$
\end{remark}

\begin{remark}\label{Rem-junction}
At this stage it is interesting to raise several questions, once we have a junction point:
\begin{itemize}
\item What is the maximum "degree" of a junction point? I.e. What is the largest number of phases that can meet at a junction point. Is this degree dimension dependent?
\item How regular is the solution $u$ to our free boundary problem close to a junction point?
\item What is the asymptotic behavior of the free boundary, close to a junction point?
\item Can two  junction points of  the same (or different) degree  come as close as possible to each other?
\end{itemize}
\end{remark}

\subsection{Multi-phase QD with triple junctions}\label{three-j}
As discussed above it is far from obvious if there are multi-phase QD with a three junction point. We shall now prove that such QD can be constructed.  We shall  confine the construction  to two-space dimension.

We start with three non-negative measures
$$
\mu_i = \chi_{B_{1/5}(z_i) }, \quad \hbox{where }
z_i=(\cos ((2i-1)\pi/3),\sin((2i-1)\pi/3)),\;\;\hbox{and}\;\; i=1,2,3.
 $$
We can solve (in a standard way) a one phase  minimization problem (obstacle type)
$$
\int_{\Pi} \frac{1}{2} |\nabla u|^2 - f_j u
$$
in the large cone
$$\Pi=\{ x: \ x= (r\cos (\theta), r\sin(\theta )) , \  0<\theta < 2\pi/3  \},$$
with force term  $f_j=- 1 + j\cdot\mu_1 $ and zero  boundary value  on $\partial \Pi$. This problem has a solution $u^j$, with bounded $\Omega^j=\{u^j >0\}$.
Next consider an even reflection of
$u^j$ in both boundaries rays of $\Pi$. This gives us two symmetric copies of $u^j$, and $\Omega^j$.

Moreover it is not hard to see that $u=(u^j_1, u^j_2,u^j_3)$ is a solution to our three-phase problem and along with $\Omega^j_i = \{u^j_i >0 \}$ gives us a three phase QD for the measures $\mu_i,$ where $i=1,2,3.$

Now let $v= x_2(x_2 + \sqrt3 x_1)/2 $ in the cone $\Pi$, and observe that $v\geq 0 $ there, and $\Delta v = 1$, along with 
zero boundary values on $\partial \Pi$.\footnote{It should be remarked that  this approach will fail for proving existence of  a four-junction point, due to the simple fact that we cannot create explicit examples like the solution v in a cone with an angle small enough to allow four junctions or more and still having $v\geq 0$ in its support.}
It is also apparent that
\begin{equation}\label{large-values}
\exists \ j_0 \ :  \quad   u^j \geq v  \hbox{ on } \partial B_{3/4} \cap \Pi
\quad \forall \ j \geq j_0 .
\end{equation}
Suppose for the moment this is true. Then by comparison principle, for the obstacle problem (see Theorem 3.3. in  \cite{Friedman-book}), we must have
$u^j \geq v$ in $B_{3/4}\cap \Pi$, for $j\geq j_0$. In particular then $u^j  (x) >0 $ for  $ x \in B_{3/4}\cap \Pi$ for large values of $j$. This in turn implies that   $B_{3/4}\cap \Pi \subset \Omega^j$, and we have thus created a junction point.

To close the argument we need to prove  the obvious statement \eqref{large-values}. We use the simple geometry of the problem. As long as $j$ is small, the QD is completely inside the cone, and it is a ball with center $z_1$. For some values of $j=j_1$ (not necessarily integer) this ball-solution will hit the boundary of the cone (at two points, 
$y^1, y^2$, symmetric with respect to the line $\{tz_1 , \ t > 0\}$) for the first time. Hence for all integer values of $j \geq j_1$ the positivity set   $\Omega^j$ of the minimizer $u^j$, 
 contains this ball. In particular, due to Hopf's boundary lemma, the solution  $u^j $ (when $j> j_1 + 1$) being strictly larger than $u^{j_1}$
has a non-vanishing gradient at the first touching points $y^1, y^2$.
 It is also apparent  (by the same way of comparison as above) that at any relative interior free boundary point  of 
 $(\partial \Omega^j \cap \partial \Pi ) \setminus \{0\}$ (which consists of two  an open line segments) we have  
 $|\nabla u^j |>0$.   This obviously means that on any curve going from one side of the boundary of $\Pi$ to another side and that lies  on $\partial B_{3/4} \cap \Pi$,  we have $u^j  (x) \geq   c_0 \hbox{dist} (x,\partial \Pi)  \geq v (x)$, for a universal $c_0$, provided $j$ is large enough. This proves \eqref{large-values}.

We summarize the above result in the following theorem.
\begin{theorem}\label{thm-three-phase}
 Quadrature Domains  with  triple  junction points, which does not touch the support of the measure,  do exist.
\end{theorem}

\begin{figure}[!htbp]
	\includegraphics[width=0.32\textwidth]{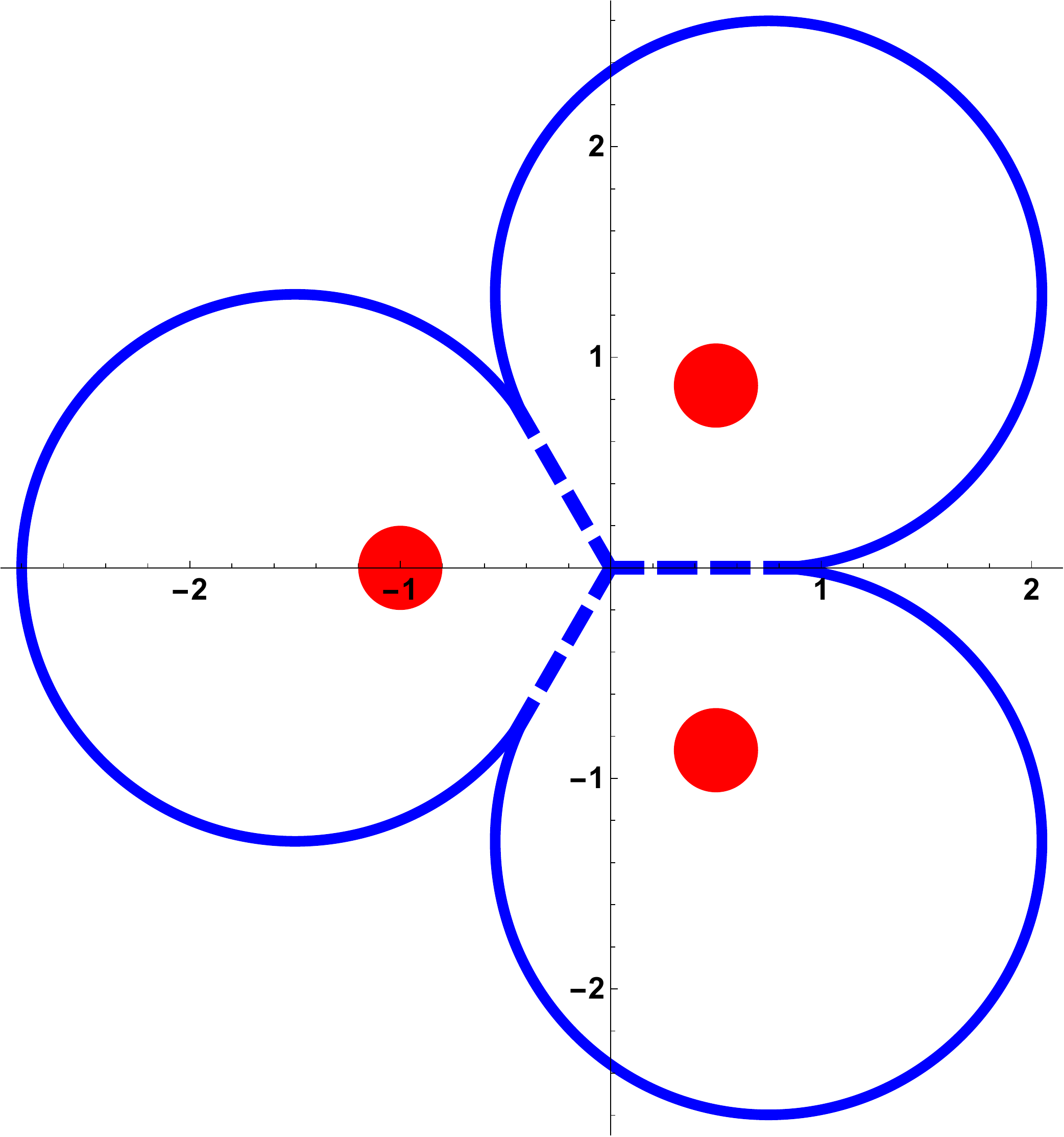}
	 \caption{\small{Three-phase QD with junction point at $(0,0).$} }
	\label{QDtriple}
\end{figure}

In Figure \ref{QDtriple} a triple junction point has been depicted. The small red disks correspond to the above defined measures $\mu_i.$ The tangential touch of  any two-phases are also shown in the picture, and this is a result of deep theorems for two-phase problems,  see \cite{PSU2012}. 

It seems plausible that one can also prove (with complete mathematical rigor) that in dimension two,  triple junction points are isolated  among points where the gradient for at least one component vanishes. This can be seen more easily at a heuristic level by a blow-up technique, and classification of global solutions. However, one needs several strong tools such as monotonicity formulas along with  complete classification of degree two homogeneous global solutions (see  below the section for Null QD),  in order to achieve such a result. This is outside the scope of the present paper.

\subsection{Multi-phase QD with  quadruple junctions}\label{sec-4-junction}
In dimension two, we shall  prove that  a multi-phase QD with four junction point does not exist, when there is a simple geometry, and provided the junction point stays away from the  support of the measure. The proof presented here does not
work in very complicated geometries, but we still believe the result should hold true.

 Let the origin be a quadruple junction point, for a given multi-phase QD, and  that for some $r>0$ we have
$B_r \cap \hbox{supp} (\mu)=\emptyset,$ where $\mu$ is the measure corresponding to that multi-phase QD. Let further  $u_i$ ($i=1,\cdots , 4$) be the corresponding function for each phase $\Omega_i$.  
Suppose we have the following simple geometry
 \begin{equation}\label{intersection}
 ( B_r \cap \partial \Omega_1 \cap \partial \Omega_3 ) \setminus \{0\} =\emptyset,
 \qquad    (B_r \cap \partial \Omega_2 \cap \partial \Omega_4 ) \setminus \{0\} =\emptyset.
 \end{equation}
   Set now $u= \lambda_3u_1 - \lambda_4u_2 + \lambda_1u_3 - \lambda_2u_4$. Then in $B_r \setminus \{0\}$, $u$ solves a two phase problem  (see Proposition  \ref{locprop} ) outside the origin, i.e. 
\begin{equation}\label{PDE-u}
\Delta u = \lambda_3\lambda_1 \chi_{\{u >0 \}} -   \lambda_4\lambda_2 \chi_{\{u <0 \}} , 
\quad \hbox{in } B_r\setminus \{0\}.
\end{equation}
It is also apparent that at the origin $\Delta u$ has no mass, so we can extend $u$ to $B_r$ and still have that it solves a two-phase problem.\footnote{This is obvious by a simple argument. Indeed  the   function $h:=v-u$, where $v$ solves   $\Delta v = f(u)$,  in $B_r$, where $f(u)$ is the right hand side in \eqref{PDE-u}, and has (some fixed bounded boundary values on $\partial B_r$), is  harmonic in $B_r\setminus \{0\}$, and bounded on $B_r$. Obviosuly this function has harmonic extension to $B_r$. }
 Hence, we can apply the theory developed for such problems (see \cite{PSU2012}). 
 Next we observe that  $\nabla u (0)=0 $, since otherwise the zero level surface  has to be a $C^{1,\alpha}$-graph, which contradicts the construction.
Hence the origin is a branch-point according to \cite{PSU2012}, and  the regularity theory for two-phase problems can be applied to conclude that both $\partial \{\pm u >0 \}$ are  locally $C^1$-graphs touching each other tangentially. This again contradicts the construction and that the  origin is a four-junction point. A similar argument can be applied to 
any junction point of order $2k$, with $k \geq 2$. It is however not clear how we can exclude junction points of order $2k + 1$ with $k>1$, i.e. odd numbers larger than three.\footnote{We remark that this argument will not work for a triple-junction point, unless at least one of the boundary curves  emanating from the origin completely separates from all other boundaries, or that gradients of all solutions on this boundary curve are zero.
This would in turn imply that the triple junction point is non-isolated among boundary points with vanishing gradient for the components. Now if this happens, then we may consider the same analysis as done in our proving non-existence of four-junction points, by grouping together two of the components (say $j=1, 3$) that are completely separated on one boundary curve and on the other they have zero gradient.  Now setting $u= \lambda_3u_1 + \lambda_1u_3 - u_2$, we shall have a solution to a two phase problem, and therefore we can use the argument done in the text above to come to a contradiction.
}

We have thus proven the following theorem.

\begin{theorem}\label{four-junction}
 Quadrature Domains  with  a quadruple  junction point, and  the geometry described in \eqref{intersection}, do  not exist.
\end{theorem}

If we allow the junction point to  hit the support of the measure, one can actually create a quadruple point as follows.

Take the following Quadrature domain $D$
\[D=\left\{(x,y)\;\mbox{s.t.}\;0<x<2,\; |y|<\frac{(1-(x-1)^2)}{2}\right\},\]
with respect to the measure $\mu$ defined as follows
\[d\mu=2(1-\sqrt{|x-1|})\chi_{[0,2]}dx.\]

\begin{figure}[!htbp]
	\includegraphics[width=0.35\textwidth]{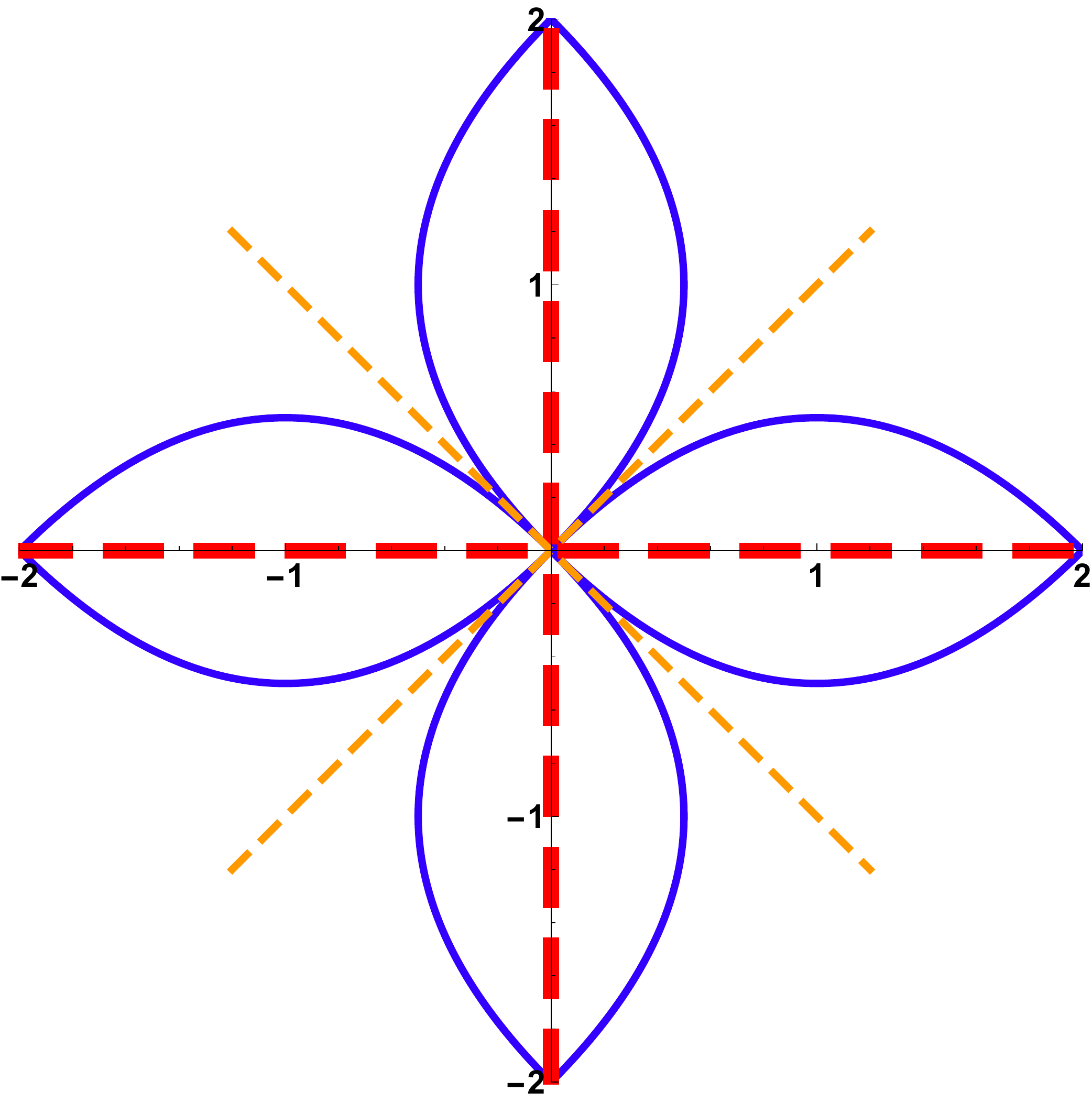}
	 \caption{\small{Quadruple  junction in the support of the measure.} }
	\label{QDfour}
\end{figure}

It can be easily verified that
\[
\int_D hdxdy=\int h d\mu, \qquad \forall h\in HL^1(D).
\]
Now, we rotate the domain $D$ clockwise with respect to the origin, by angle $\pi/2.$  Repeating this process three times we will end up with a picture as in Figure \ref{QDfour}. The orange dash-lines are the tangents at a point $(0,0).$ The red dash-lines indicate  the  measures support of corresponding rotated quadrature domain $D.$ It is clear that this new domain is a four-phase QD with quadruple junction point at $(0,0)$, but with support of the measure meeting the origin.

\subsection{Null QD with triple junctions}\label{three-j}
We shall now give explicit examples of a global three-phase  QD, with zero measure, i.e. a Null QD of three phases. 
The explicit form of the three phase version example reads as follows
\begin{align}
\begin{cases}
u_1(x,y)=\frac{\lambda_1}{2}\cdot y(y-ax) \quad & \mbox{in} \;\; {y\geq0,\; y\geq ax},\\
u_2(x,y)=-\frac{\lambda_2}{2}\cdot y(bx-y) \quad & \mbox{in} \;\; {y\leq0,\; y\leq bx},\\
u_3(x,y)=\frac{\lambda_1\lambda_2}{2(\lambda_1+\lambda_2)}\cdot (ax-y)(y-bx) \quad & \mbox{in} \;\; {bx\leq y\leq ax},
\end{cases}
\end{align}
where
\[ a=-\frac{\sqrt{\lambda_1\lambda_2+\lambda_1\lambda_3+\lambda_2\lambda_3}}{\lambda_1}, \]
and
\[ b=\frac{\sqrt{\lambda_1\lambda_2+\lambda_1\lambda_3+\lambda_2\lambda_3}}{\lambda_2}.\]
Observe that the pairing $u_{ij}:=u_i - u_j$, for $i\neq j$ satisfies a two-phase equation 
$$
\Delta u_{ij} = \lambda_i \chi_{\{u_{ij} > 0\}}- \lambda_j \chi_{\{u_{ij} < 0\}} \qquad \hbox{in } interior(\overline{\Omega_i \cup \Omega_j}),
$$
with $u_{ij}= 0 $ on $ \partial \Omega_i \cup \partial \Omega_j$, where $\Omega_i =\{u_i >0\}$.

A classification of global solutions (Null QD) seems to be a key element in any local analysis of junction points, as discussed at the beginning of this section, see  Remark \ref{Rem-junction}. For the moment, lack of techniques to analyze such problems prevent us from any further study of the junction points.

\section{An Application from control theory}
We shall now  present an application of our problem which can be related to  optimal control theory.
 Let $ U_{ad} $ be the admissible class consisting of all vectors $h=(h_1,\cdots , h_m)$, according to
$$U_{ad}:= \left\{ h_i  \in L^\infty (\Omega): \   0 \leq   h_i \leq \lambda_i , \  h_i h_j=0 , \hbox{ for } i\neq j, \ \{ h_i >0 \} \supset \hbox{supp} (\mu_i)  \right\} , $$
where $\mu_i$ are given measures, with mutually disjoint supports.
 Let  further  $\Omega_i :=\{ h_i >0 \} $, and $v_i$ be
the solution of the linear problem
\begin{equation}\label{optimal-control-application}
\left\{
\begin{aligned}
&\Delta v_i =h_i   -  \mu_i  \quad \text{in } \Omega_i  \\
&v_i=0  \quad \text{on}\  \partial \Omega_i,
\end{aligned}
\right.
\end{equation}
and extended as zero outside of $\Omega_i.$
Here  $h_i$ is a control function.
Define
$$
I(h_1,\cdots , h_m):=  \sum_{i=1}^m \int\limits_{\Omega }|\nabla v_i|^2+ v_i^+ \lambda_i - v_i\mu_i  ,
$$
for all $h_i  \in
U_{ad}$. Observe that the functional depends on $(h_1, \cdots , h_m)$ implicitly through $v_i$.
\newline
Using integration by parts for the gradient part we easily calculate that

\begin{align*}
I(h_1,\cdots , h_m)&=\sum_{i=1}^m\int\limits_{\Omega }\left(-\Delta v_i\cdot v_i + v_i^+ \lambda_i - v_i\mu_i\right)\\
&=\sum_{i=1}^m\int\limits_{\Omega_i }\left( (\mu_i-h_i )v_i+v_i^+ \lambda_i - v_i\mu_i\right)\\&=
\sum_{i=1}^m\int\limits_{\Omega_i }\left(v_i^+(\lambda_i - h_i)+h_iv_i^-\right)\geq 0.
\end{align*}
Here $I(h)=0$ iff $h_i =\lambda_i \chi_{ \{v_i >0\}}$, and $v_i^-=0$. Hence  $h_i = \lambda_i \chi_{\{v_i >0\}} $, with $v_i \geq 0$, minimizes the functional $I$ if $(v_1, \cdots , v_m)$ is a minimizer  to our functional.

\bibliographystyle{acm}

\bibliography{MultiQD}

\end{document}